\documentclass[MSNbibl,number,citesort,seceqn,dvips]{arxbj}
\usepackage{graphicx}

\aid{0}
\volume{18}
\issue{2}
\pubyear{2012}
\firstpage{735}
\lastpage{745}
\doi{10.3150/11-BEJ352}

\makeatletter
\DeclareMathAlphabet{\mathbbl}  {U}{mt2hrb}{m}{n}
\SetMathAlphabet{\mathbbl}{bold}{U}{mt2hrb}{b}{n}

\renewcommand{\epsilon}{\varepsilon}

\newtheorem{theorem}{Theorem}[section]
\newremark{ex}{Example}[section]
\makeatother

\begin{document}
\begin{frontmatter}

\title{On moving-average models with feedback}
\runtitle{Moving-average models}

\begin{aug}
%%%% inicialai - be tarpu
\author[1]{\fnms{Dong} \snm{Li}\corref{}\thanksref{1,e1}\ead[label=e1,mark]{malidong@ust.hk}},
\author[1]{\fnms{Shiqing} \snm{Ling}\thanksref{1,e2}\ead[label=e2,mark]{maling@ust.hk}}
\and
\author[2]{\fnms{Howell} \snm{Tong}\thanksref{2}\ead[label=e3]{howell.tong@gmail.com}}

\runauthor{D. Li, S. Ling and H. Tong}
\address[1]{Department of Mathematics, Hong Kong University of Science and Technology, Hong
Kong.\\
\printead{e1,e2}}
\address[2]{Department of Statistics, London School of Economics and Political
Science, Houghton Street, London, WC2A 2AE, UK. \printead{e3}}
\end{aug}

% HISTORY:
\received{\smonth{7} \syear{2010}}
\revised{\smonth{12} \syear{2010}}

% ABSTRACT
\begin{abstract}
Moving average models, linear or nonlinear, are characterized by
their short memory. This paper shows that, in the presence of
feedback in the dynamics, the above characteristic can disappear.
\end{abstract}

% KEYWORDS
\begin{keyword}
\kwd{ACF}
\kwd{ergodicity}
\kwd{existence}
\kwd{feedback}
\kwd{leptokurticity}
\kwd{memory}
\kwd{stationarity}
\kwd{thresholds}
\end{keyword}

\end{frontmatter}

%s1 #&#
\section{Introduction}

Since the introduction by Slutsky \cite{ref15}, moving average models have
played a significant role in time series analysis, especially in
finance and economics. The models have been extended to include
measurable (nonlinear) functions of independent and identically
distributed random variables, representing unobservable and purely
random impulses, for example, Robinson \cite{ref14}. The characterizing feature of
these models is the cut-off of the auto-covariance functions when
they exist, implying that they are models of short memory. Another
interesting feature of these models is the homogeneity of the random
impulses, free from any feedback in the generating mechanism. Now,
Slutsky developed these models in an economic context; the random
impulses may correspond to unobservable political factors. In such a
context, as well as in other contexts for which these models are
relevant (e.g., business studies), it can be argued that feedback is
often present: political decisions are often predicated on economic
conditions. One simple way to incorporate feedback in these models
is through the notion of thresholds, that is, on-off feedback
controllers.

Since Tong \cite{ref16} initiated the threshold notion in time series
modelling, the notion has been extensively used in the literature,
especially for the threshold autoregressive (TAR) or TAR-type
models. For these models, some basic and probabilistic properties
were given in Chan \textit{et al.} \cite{ref3} and Chan and Tong \cite{ref4}.
More related results can be found in  An and Huang \cite{ref1}, Brockwell
\textit{et al.} \cite{ref2}, Chen and Tsay \cite{ref5}, Cline and Pu \cite{ref6,ref7}, Ling \cite{ref9},
Ling \textit{et al.} \cite{ref11}, Liu and Susko \cite{ref12}
and Lu \cite{ref13}, among others. A fairly comprehensive review of
threshold models is available in Tong \cite{ref17} and a selective survey
of the history of threshold models is given by Tong
\cite{ref18}.\vadjust{\goodbreak}

However, most work to-date on the threshold model  has primarily
concentrated on the TAR or the TAR-type model. The threshold moving
average (TMA) model, that is a moving average model with a simple
on-off feedback control mechanism, has not attracted as much
attention. As far as we know, only a few results are available for
the TMA model. Brockwell \textit{et al.} \cite{ref2} investigated a
threshold autoregressive and moving-average (TARMA) model and
obtained a~strictly stationary and ergodic solution to the model
when the MA part does not contain any threshold component.
Unfortunately, their TARMA model does not cover the TMA model as a
special case. Using the Markov chain theory, Liu and Susko \cite{ref12}
provided the existence of the strictly stationary solution to the
TMA model without any restriction on the coefficients. However, they
neither gave an explicit form of the solution nor proved the
ergodicity.  A similar result can be found in Ling \cite{ref9}. Ling
\textit{et al.} \cite{ref11} gave a sufficient condition for the ergodicity of the
solution for a first order TMA model under some restrictive
conditions. These results have been extended to the first-order TMA
model with more than two regimes.  However, the uniqueness and the
ergodicity of the solution are still open problems for
higher-order TMA models.

In this paper, we use a different approach to study the TMA model
without resorting to the Markov chain theory.  Note that the TMA
model involves a feedback control mechanism. An intuitive and simple
idea is to seek a closed form of the solution in terms of the above
mechanism, which is expressible as an indicator function. We can
show that for the TMA model there always exists a unique strictly
stationary and ergodic solution without any restriction on the
coefficients of the TMA model. More importantly, for the first time
in the literature, an explicit/closed form of the solution is
derived. In addition, for the correlation structure, we show that
the ACF (when it exists) of the TMA model typically does not cut
off. In fact, it has a  much richer structure. For example, it can
exhibit almost long memory, although it generally decays at an
exponential rate. Furthermore, the difference between the joint
two-dimensional distribution and the corresponding product of its
marginal distributions also decays to zero at an exponential rate as
the lag tends to infinity.

The rest of the paper is organized as follows. Section \ref{sec2} discusses
the strict stationarity and ergodicity of the TMA model. Section \ref{sec3}
studies the asymptotic behaviour of the ACF of the TMA model and
other correlation structure. We conclude in Section \ref{sec4}. All proofs of
the theorems are relegated to the \hyperref[appendix]{Appendix}.
%%%%%%%%%%%%%%%%%%%%%%%%%%%%%%%%%%%%%%%%%%%%%%%%%%%%%%%%%%%%%%%%%%%%%%%%%%%%%%%%%%%%%%%%%%%%%%%%%%%%%%%%%%%%%%%%%%%%%%%%%%
%                                                                                                                      %%%
%             Section 2: Stationarity and ergodicity of $\operatorname{TMA}(q)$ models                                                %%%
%                                                                                                                      %%%
%%%%%%%%%%%%%%%%%%%%%%%%%%%%%%%%%%%%%%%%%%%%%%%%%%%%%%%%%%%%%%%%%%%%%%%%%%%%%%%%%%%%%%%%%%%%%%%%%%%%%%%%%%%%%%%%%%%%%%%%%%
%s2 #&#
\section{Stationarity and ergodicity of $\operatorname{TMA}(q)$
models}\label{sec2}

We first  consider  a $\operatorname{TMA}(q)$ model which satisfies the following
equation:
%
%e1 #&#
\begin{eqnarray}\label{p}
y_n= \cases{\displaystyle
\mu_1+e_n+\sum_{i=1}^q\phi_ie_{n-i},& \quad if
$y_{n-d}\leq r$,\vspace*{2pt}\cr
\displaystyle\mu_2+e_n+\sum_{i=1}^q\psi_ie_{n-i},& \quad if $y_{n-d}>r$,
}
\end{eqnarray}
where $\{e_n\}$ is a sequence of i.i.d. random variables. Here, $q$
and $ d $ are positive integers,  $r\in\mathbb{R}$, the real line,
is the threshold parameter, and $\mu_1, \mu_2, \phi_i$ and $\psi_i$,
$i=1,\ldots,q$,  are real coefficients.

For the sake of simplicity, we adopt the following notation:
\begin{eqnarray*}
U_n=\mathbbl{1}(a_n\leq r)\quad\mbox{and}\quad
W_n=\mathbbl{1}(b_n\leq r)-\mathbbl{1}(a_n\leq
r),
\end{eqnarray*}
where $\mathbbl{1}(\cdot)$ is an indicator function,
%
%e2 #&#
\begin{eqnarray}\label{op}
a_n=\mu_2+e_n+\sum_{i=1}^q\psi_ie_{n-i}\quad\mbox{and}\quad
b_n=\mu_1+e_n+\sum_{i=1}^q\phi_ie_{n-i}.
\end{eqnarray}
The following  theorem gives  the strict stationarity and ergodicity
of model (\ref{p}).\vspace*{-3pt}

%%%%%%%%%%%%%%%%%%%%%%%%%%%%%%%%%%%%%%%%%%%%%%%%%%%%%%%%%%%%%%%%%%%%%%%%%%%%%%%%%%%%%%%%%%%%%%%%%%%%%%%%%%%%%%%%%%%%%%%%%%
%                                                                                                                      %%%
%                                        Theorem 2.1                                                                   %%%
%                                                                                                                      %%%
%%%%%%%%%%%%%%%%%%%%%%%%%%%%%%%%%%%%%%%%%%%%%%%%%%%%%%%%%%%%%%%%%%%%%%%%%%%%%%%%%%%%%%%%%%%%%%%%%%%%%%%%%%%%%%%%%%%%%%%%%%
%t1 #&#
\begin{theorem}\label{sethm1}
Suppose that $\{e_n\}$ is a sequence of i.i.d. random variables with
$\mathbb{P}(a_n\leq r,\allowbreak b_n\leq r)+\mathbb{P}(a_n> r, b_n>r)\not=0$.
Then $y_n$ has a unique strictly stationary and ergodic solution
expressed by
\[
y_n=\mu_2+e_n+\sum_{i=1}^q\psi_ie_{n-i}+\Biggl[(\mu_1-\mu_2)+\sum_{i=1}^q(\phi_i-\psi_i)e_{n-i}\Biggr]\alpha_{n-d},\qquad
\mbox{a.s.},
\]
where
\[
\alpha_{n-d}=\sum_{j=1}^{\infty}\Biggl[\Biggl(\prod_{s=1}^{j-1}W_{n-sd}\Biggr)U_{n-jd}\Biggr],
\qquad \mbox{in $L^1$ and a.s.}\vspace*{-3pt}
\]
\end{theorem}
If $e_{1}$ has  a strictly and continuously positive density on
$\mathbb{R}$ (e.g.,  normal, Student's $t_v$ or double
exponential distribution), then $\mathbb{P}(a_n\leq r, b_n\leq
r)+\mathbb{P}(a_n> r, b_n>r)\not=0$. The basic  idea for Theorem
\ref{sethm1} is a direct and concrete expression in terms of
$\mathbbl{1}(y_{n-d}\leq r)$, without resorting to the Markov chain
theory. Theorem \ref{sethm1} shows that the $\operatorname{TMA}(q)$
model is always
stationary and ergodic as is the $\operatorname{MA}(q)$ model.\vspace*{-3pt}

%s3 #&#
\section{The ACF of $\operatorname{TMA}(q)$ models}\label{sec3}\vspace*{-3pt}

The ACF plays a crucial role in studying the correlation structure
of weakly stationary time series. It is well known that for a causal
$\operatorname{ARMA}(p, q)$ model, its ACF $\rho_k$ goes to zero at an exponential
rate as $k$ diverges to infinity. The exact formula for ACF can be
obtained although its closed form is not compact. However, for a
general nonlinear time series model, it is rather difficult to
obtain an exact formula for the ACF and to study the asymptotic
behaviour. Additionally, the notion of memory, short or long, is
closely associated with the ACF. One significant fact is that a
causal $\operatorname{ARMA}(p, q)$ model is short-memory. For a general nonlinear
time series model, due to its complicated structure, there is no
universally accepted criterion for determining whether or not it is
short-memory. As for some specific time series model, an ad hoc
approach is usually adopted.

One important characteristic of the $\operatorname{MA}(q)$ model is that its ACF
cuts off after lag $q$. Interestingly, this property is not\vadjust{\goodbreak}
generally inherited by the TMA model; this is not surprising
theoretically because the TMA model involves some nonlinear
feedback. Another interesting fact is that although a TMA model is
generally short-memory, in some cases it can exhibit some almost
long-memory phenomena; see Example \ref{seex4}. The following
theorem characterizes the ACF of model (\ref{p}).\vspace*{-3pt}

%%%%%%%%%%%%%%%%%%%%%%%%%%%%%%%%%%%%%%%%%%%%%%%%%%%%%%%%%%%%%%%%%%%%%%%%%%%%%%%%%%%%%%%%%%%%%%%%%%%%%%%%%%%%%%%%%%%%%%%%%%
%                                                                                                                      %%%
%                                        Theorem 3.1                                                                   %%%
%                                                                                                                      %%%
%%%%%%%%%%%%%%%%%%%%%%%%%%%%%%%%%%%%%%%%%%%%%%%%%%%%%%%%%%%%%%%%%%%%%%%%%%%%%%%%%%%%%%%%%%%%%%%%%%%%%%%%%%%%%%%%%%%%%%%%%%
%t2 #&#
\begin{theorem}\label{sethm2}
Suppose that the condition in Theorem \textup{\ref{sethm1}} is satisfied and
$\mathbb{E}|e_1|^{2}<\infty$.  Then there exists a constant $\rho\in
(0,1)$ such that $\rho_k=\mathrm{O}(\rho^{k})$.\vspace*{-3pt}
\end{theorem}

Theorem \ref{sethm2} indicates that the TMA model (\ref{p}) is
short-memory. The next theorem describes the relationship between
the two-dimensional joint distribution and the corresponding marginal
distributions.\vspace*{-3pt}
%%%%%%%%%%%%%%%%%%%%%%%%%%%%%%%%%%%%%%%%%%%%%%%%%%%%%%%%%%%%%%%%%%%%%%%%%%%%%%%%%%%%%%%%%%%%%%%%%%%%%%%%%%%%%%%%%%%%%%%%%%
%                                                                                                                      %%%
%                                        Theorem 3.2                                                                   %%%
%                                                                                                                      %%%
%%%%%%%%%%%%%%%%%%%%%%%%%%%%%%%%%%%%%%%%%%%%%%%%%%%%%%%%%%%%%%%%%%%%%%%%%%%%%%%%%%%%%%%%%%%%%%%%%%%%%%%%%%%%%%%%%%%%%%%%%%
%t3 #&#
\begin{theorem}\label{sethm3}
Suppose that $\{e_n\}$ is i.i.d. random variables having a
continuously, boundedly and strictly positive density. Then, for any
$u, v\in\mathbb{R}$ and $k\geq 1$, there exists a constant $\rho\in
(0,1)$ such that
\[
|\mathbb{P}(y_0\leq u, y_k\leq v)-\mathbb{P}(y_0\leq
u)\mathbb{P}(y_k\leq v)|=\mathrm{O}(\rho^{k}).\vspace*{-3pt}
\]
\end{theorem}

Actually, Theorem \ref{sethm3} still holds for
$\operatorname{Cov}(\mathbbl{1}(u_1<y_0\leq u_2), \mathbbl{1}(v_1<y_k\leq v_2))$
where $-\infty\leq u_1<u_2\leq \infty$ and $-\infty\leq v_1<v_2\leq
\infty$.
%%%%%%%%%%%%%%%%%%%%%%%%%%%%%%%%%%%%%%%%%%%%%%%%%%%%%%%%%%%%%%%%%%%%%%%%%%%%%%%%%%%%%%%%%%%%%%%%%%%%%%%%%%%%%%%%%%%%%%%%%%
%                                                                                                                      %%%
%                                          Section 4: Some examples                                                    %%%
%                                                                                                                      %%%
%%%%%%%%%%%%%%%%%%%%%%%%%%%%%%%%%%%%%%%%%%%%%%%%%%%%%%%%%%%%%%%%%%%%%%%%%%%%%%%%%%%%%%%%%%%%%%%%%%%%%%%%%%%%%%%%%%%%%%%%%%
Next, we consider some special TMA models and  study their ACFs as
well as some other properties.\vspace*{-3pt}
%%%%%%%%%%%%%%%%%%%%%%%%%%%%%%%%%%%%%%%%%%%%%%%%%%%%%%%%%%%%%%%%%%%%%%%%%%%%%%%%%%%%%%%%%%%%%%%%%%%%%%%%%%%%%%%%%%%%%%%%%%
%                                                                                                                      %%%
%                                        Example 4.1                                                                   %%%
%                                                                                                                      %%%
%%%%%%%%%%%%%%%%%%%%%%%%%%%%%%%%%%%%%%%%%%%%%%%%%%%%%%%%%%%%%%%%%%%%%%%%%%%%%%%%%%%%%%%%%%%%%%%%%%%%%%%%%%%%%%%%%%%%%%%%%%
\begin{ex}\label{seex1}
Suppose that $y_n$ is defined as
\[
y_n= \cases{
\mu_1+e_n,& \quad  if $y_{n-1}\leq r$,\cr
\mu_2+e_n, & \quad if $y_{n-1}>r$,
}
\]
where $\{e_n\}$ satisfies the condition in Theorem \ref{sethm1} with
mean $0$ and finite variance $\sigma^2$.\vspace*{-3pt}
\end{ex}

This example can also be regarded as a special case of the TAR
model, which was studied in Tong \cite{ref17},  Question 29, page 212. By
calculation, we have the ACF of $\{y_n\}$
\[
\rho_k=\frac{(\mu_1-\mu_2)\lambda_k+(\mu_1-\mu_2)^2\delta_0(1-\delta_0)\beta^k}{\sigma^2+(\mu_1-\mu_2)^2\delta_0(1-\delta_0)}
\qquad \mbox{for } k\geq 1,
\]
where $\lambda_k=\mathbb{E}[e_{n-k}\mathbbl{1}(y_{n-1}\leq r)]$,
$\beta=[G(r-\mu_1)-G(r-\mu_2)]\in (-1, 1)$ and
$\delta_0=G(r-\mu_2)/\allowbreak[1-G(r-\mu_1)+G(r-\mu_2)]$. Here, $G(x)$ is the
distribution function of $e_1$.

\begin{figure}

\includegraphics{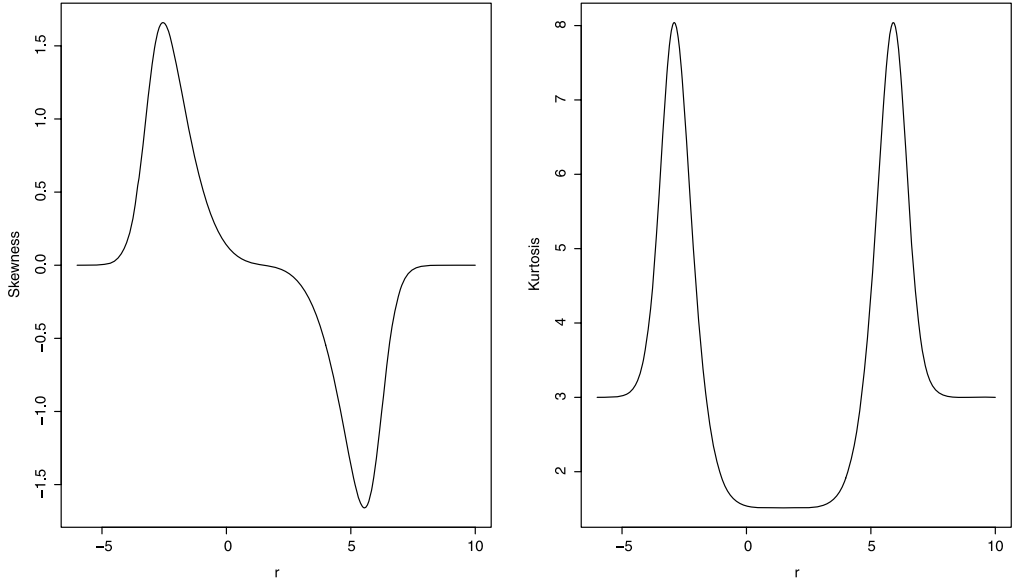}

\caption{The skewness (left) and the kurtosis
(right) of $y_n$ as functions of $r$ when $e_1$ is standard
normal.}\label{sefig0}
\end{figure}

Clearly, the ACF does not possess the cut-off property except for
$\mu_1=\mu_2$. Generally, $\rho_k$ decays exponentially since
$\lambda_k=\mathrm{O}(\rho^k)$ for some $\rho\in(0, 1)$ by the proof of
Theorem \ref{sethm3}. In the nonlinear time series literature, the
search for a nonlinear AR model with long memory has been largely in
vain. Against this background, it is interesting to note that as
$\mu_1\rightarrow\infty$ and $\mu_2\rightarrow-\infty$, $\rho_k$ can
exhibit almost long memory in that $\rho_k$ can be made to decay
arbitrarily slowly. Note that Example \ref{seex1} can be driven by a
white noise process with a thin tailed distribution. The skewness
and the kurtosis of $y_n$ are also available explicitly and
interesting. Specifically,
\[
\mathrm{skewness}=\frac{\mathbb{E}e_1^3+(\mu_1-\mu_2)^3(\delta_0-3\delta_0^2+2\delta_0^3)}
{[\sigma^2+(\mu_1-\mu_2)^2\delta_0(1-\delta_0)]^{3/2}}\vadjust{\goodbreak}
\]
and
\[
\mathrm{kurtosis}=\frac{\mathbb{E}e_1^4+6\sigma^2(\mu_1-\mu_2)^2\delta_0(1-\delta_0)
+(\mu_1-\mu_2)^4(\delta_0-4\delta_0^2+6\delta_0^3-3\delta_0^4)}
{[\sigma^2+(\mu_1-\mu_2)^2\delta_0(1-\delta_0)]^{2}},
\]
respectively. The impact of the threshold parameter $r$ is related
to the bi-modality of the marginal density, which can be established
by simple calculation. When $e_1$ is standard normal and
$(\mu_1,\mu_2)=(4,-1)$, Figure \ref{sefig0} shows the skewness and the
kurtosis of $y_n$ as functions of $r$.

%f1 #&#

\begin{ex}\label{seex2}
Suppose that $\{y_n\}$ follows a $\operatorname{TMA}(1)$ model without drift:
\[
y_n= \cases{
e_n+\phi e_{n-1},& \quad if  $y_{n-2}\leq r$,\cr
e_n+\psi e_{n-1},& \quad if $y_{n-2}>r$,
}
\]
where $\{e_n\}$ satisfies the condition in Theorem \ref{sethm1},
having zero mean and finite variance.
\end{ex}

After simple calculation, we have the ACF
\[
\rho_k= \cases{
\displaystyle\frac{\psi+(\phi-\psi)\varrho}{1+\psi^2+(\phi^2-\psi^2)\varrho},&\quad if $k=1$,\vspace*{2pt}\cr
0,& \quad if $k\geq 1$,
}
\]
where $\varrho=\mathbb{P}(e_2+\psi e_{1}\leq r)/[\mathbb{P}(e_2+\phi
e_{1}> r)+\mathbb{P}(e_2+\psi e_{1}\leq r)]\in[0, 1)$.

This example shows that for some special $\operatorname{TMA}(q)$ model, the ACF may
be cut off after lag $q$. In particular, if $\phi=\psi$, then the
ACF coincides with that of the classical linear MA(1) model.
Unfortunately, for general TMA models with $d\leq q$, there are no
explicit expressions available for the ACFs due to the extremely
complicated dependence of $y_{t-d}$ on $\{e_{t-j}, d\leq j\leq q\}$.
However, we can obtain the sample ACFs of TMA models by simulation.

%%%%%%%%%%%%%%%%%%%%%%%%%%%%%%%%%%%%%%%%%%%%%%%%%%%%%%%%%%%%%%%%%%%%%%%%%%%%%%%%%%%%%%%%%%%%%%%%%%%%%%%%%%%%%%%%%%%%%%%%%%
%                                                                                                                      %%%
%                                        Example 4.4                                                                   %%%
%                                                                                                                      %%%
%%%%%%%%%%%%%%%%%%%%%%%%%%%%%%%%%%%%%%%%%%%%%%%%%%%%%%%%%%%%%%%%%%%%%%%%%%%%%%%%%%%%%%%%%%%%%%%%%%%%%%%%%%%%%%%%%%%%%%%%%%
\begin{ex}\label{seex4}
 Suppose that $\{y_n\}$ follows the $\operatorname{TMA}(1)$ model:
%e3 #&#
\begin{equation}\label{llt}
y_n= \cases{
5+e_n+0.2 e_{n-1},& \quad if $y_{n-1}\leq 0.5$,\cr
-3+e_n+0.8 e_{n-1},& \quad if $y_{n-1}>0.5$,
}
\end{equation}
where $\{e_n\}$ is  i.i.d. standard normal.
\end{ex}

\begin{figure}

\includegraphics{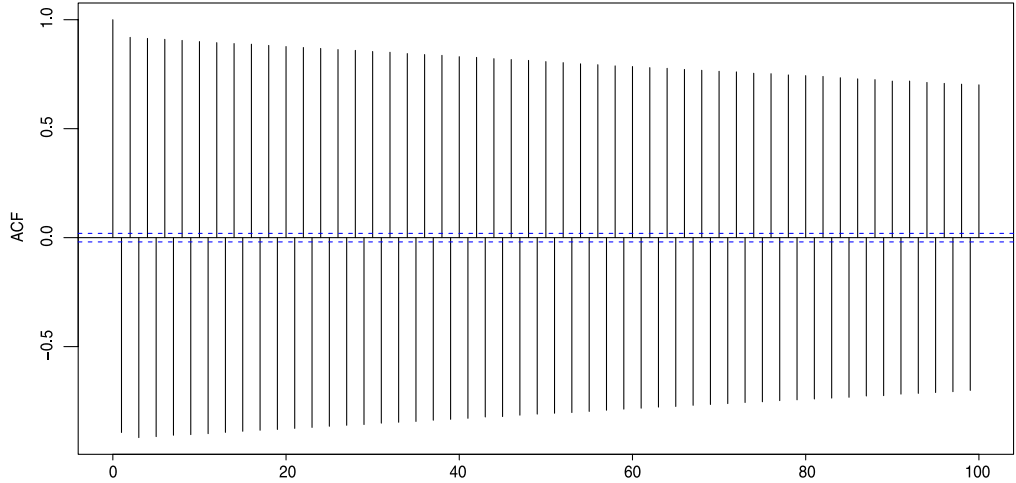}%

\caption{The sample  ACFs of model (\protect\ref{llt}).}\label{sefig2}

\end{figure}

This model produces a time series that mimics a unit root and long
memory. In Figure~\ref{sefig2}, the sample ACF of model (\ref{llt})
decays slowly, although model (\ref{llt}) is stationary.

%f2 #&#

%s4 #&#
\section{Concluding remarks}\label{sec4}

Conventional moving average models, whether linear or nonlinear,
assume absence of any feedback control mechanism. This paper shows
that the introduction of simple feedback can enrich the structure of
moving average models. For example, their ACF need not cut off but
can now exhibit (near) long memory. Their distributions can be
leptokurtic even when driven by Gaussian white noise. In nonlinear
time series modeling, moving average models have been overshadowed
by autoregressive models. Our study suggests that, by introducing a
simple feedback mechanism, the notion of moving average possesses
some unexpected properties beyond the shadow.

\begin{appendix}\label{appendix}
%s5 #&#
\section{Proofs of theorems}\label{app}
%%%%%%%%%%%%%%%%%%%%%%%%%%%%%%%%%%%%%%%%%%%%%%%%%%%%%%%%%%%%%%%%%%%%%%%%%%%%%%%%%%%%%%%%%%%%%%%%%%%%%%%%%%%%%%%%%%%%%%%%%%
%                                                                                                                      %%%
%                                Appendix A.1: Proofs of Theorem 2.1                                                   %%%
%                                                                                                                      %%%
%%%%%%%%%%%%%%%%%%%%%%%%%%%%%%%%%%%%%%%%%%%%%%%%%%%%%%%%%%%%%%%%%%%%%%%%%%%%%%%%%%%%%%%%%%%%%%%%%%%%%%%%%%%%%%%%%%%%%%%%%%
%s5.1 #&#
\subsection{\texorpdfstring{Proof of Theorem \protect\ref{sethm1}}{Proof of Theorem 2.1}}
From model $(\ref{p})$, $\mathbbl{1}(y_{n}\leq
r)=U_n+W_n\mathbbl{1}(y_{n-d}\leq r).$  Iterating $k\ge 1$ steps, we
have
\[
\mathbbl{1}(y_{n}\leq
r)=\sum_{j=0}^{k-1}\Biggl[\Biggl(\prod_{s=0}^{j-1}W_{n-sd}\Biggr)U_{n-jd}\Biggr]+\Biggl(\prod_{i=0}^{k-1}W_{n-id}\Biggr)\mathbbl{1}(y_{n-kd}\leq
r)
\]
with the convention $\prod_{0}^{-1}=1$.  Let
\[
\alpha_{n,k}=\sum_{j=0}^{k-1}\Biggl[\Biggl(\prod_{s=0}^{j-1}W_{n-sd}\Biggr)U_{n-jd}\Biggr].
\]
For given $d$ and $q$, there exists a unique nonnegative integer $m$
such that $md<\max(d, q+1)\leq (m+1)d$. Let
$\delta=\mathbb{E}|W_{1}|$. Under the condition in Theorem
\ref{sethm1}, it is not difficult to prove that $0\leq\delta<1$.
Observing that both $\{U_n\}$ and $\{W_n\}$ are $q$-dependent
sequences,  we can extract an independent subsequence
$\{W_{n-j(m+1)d}, j=0, 1,\ldots, [\frac{k-1}{m+1}]\}$ from the
sequence $\{W_{n-id}, i=0, 1, 2,\ldots,k-1\}$, where $[a]$ denotes the
integral part of $a$.  Since $|U_n|\le 1$ and $|W_n|\le 1$,  it
yields that
\[
\mathbb{E}\Biggl|\Biggl(\prod_{i=0}^{k-1}W_{n-id}\Biggr)U_{n-kd}\biggr|\leq(\mathbb{E}|W_{1}|)^{[(k-1)/(m+1)]},
\]
implying
\[
\sum_{j=1}^\infty\mathbb{E}\Biggl|\Biggl(\prod_{i=0}^{j-1}W_{n-id}\Biggr)U_{n-jd}\Biggr|\leq
\sum_{j=1}^\infty\delta^{[(j-1)/(m+1)]}=(m+1)\sum_{k=0}^\infty\delta^k<\infty.
\]
Using the above inequalities, we can prove that $
\mathbb{E}|\alpha_{n,s}-\alpha_{n,t}|\rightarrow0$ as $s,
t\rightarrow\infty$ for each fixed $n$. By the Cauchy criterion,
$\alpha_{n,k}$ converges in $L^1$ as $k\rightarrow\infty$. Write the
limit as
\[
\alpha_{n}=\sum_{j=0}^{\infty}\Biggl[\Biggl(\prod_{s=0}^{j-1}W_{n-sd}\Biggr)U_{n-jd}\Biggr].
\]
Applying the inequalities  above again, it is easy to get
\[
\sum_{k=1}^\infty\mathbb{E}|\alpha_{n,k}-\alpha_{n}|
\leq\sum_{k=1}^\infty\sum_{j=k}^\infty\delta^{[(j-1)/(m+1)]}<\infty,
\]
yielding that
\[
\lim_{k\rightarrow\infty}\alpha_{n, k}=\alpha_n, \qquad \mbox{in
$L^1$ and a.s.}
\]

Furthermore, recall that $U_n=\mathbbl{1}(a_n\leq r)$ and
$W_n=\mathbbl{1}(b_n\leq r)-\mathbbl{1}(a_n\leq
r)$, where $a_n$ and~$b_n$ are defined in (\ref{op}),  we have
the iterative sequence: $\alpha_{n,1}=U_n$ and
\[
\alpha_{n,k}=U_n+W_n\alpha_{n-d,k-1}=(1-\alpha_{n-d,k-1})\mathbbl{1}(a_n\leq
r)+\alpha_{n-d,k-1}\mathbbl{1}(b_n\leq r)
\]
for each $n$ and $k\geq 1$. Note that $\alpha_{n,k}$ and
$\alpha_{n-d,k}$ have the same distribution for fixed~$k$ since the
error $\{e_i\}$ is i.i.d. By induction over $k$, we have that
$\alpha_{n,k}$ only takes two values~0 and 1 a.s. since
$\alpha_{n,1}$ only takes 0 and 1. Thus, $\alpha_n$ at most takes
two values 0 and 1 a.s., namely, $\alpha_n=\mathbbl{1}(\alpha_n=1)$
a.s. Define a new sequence $\{S_n\}$
\[
S_n=\mu_2+e_n+\sum_{i=1}^q\psi_ie_{n-i}+\Biggl[(\mu_1-\mu_2)+\sum_{i=1}^q(\phi_i-\psi_i)e_{n-i}\Biggr]\alpha_{n-d}.
\]
By simple calculation, we have
\begin{eqnarray*}
\mathbbl{1}(S_n\leq r)&=&\mathbbl{1}(a_n\leq
r)\mathbbl{1}(\alpha_{n-d}=0)+\mathbbl{1}(b_n\leq
r)\mathbbl{1}(\alpha_{n-d}=1)\\
&=&U_n+W_n\mathbbl{1}(\alpha_{n-d}=1)\\
&=&U_n+W_n\alpha_{n-d}=\alpha_n,\qquad \mbox{a.s.}
\end{eqnarray*}
Hence,
\[
S_n=\mu_2+e_n+\sum_{i=1}^q\psi_ie_{n-i}+\Biggl[(\mu_1-\mu_2)+\sum_{i=1}^q(\phi_i-\psi_i)e_{n-i}\Biggr]\mathbbl{1}(S_{n-d}\leq
r),\qquad \mbox{a.s.}
\]
Thus, $\{S_{n}\}$ is  the solution of model (2.1) which is strictly
stationary and ergodic.

To uniqueness, suppose that $\tilde{S}_n$ is a solution to model
(\ref{p}), then
\[
\mathbbl{1}(\tilde{S}_{n}\leq
r)=U_n+W_n\mathbbl{1}(\tilde{S}_{n-d}\leq r).
\]
Iterating the above equation, one can get for $k\geq 1$
\[
\mathbbl{1}(\tilde{S}_{n}\leq
r)=\alpha_{n,k}+\Biggl(\prod_{i=0}^{k-1}W_{n-id}\Biggr)\mathbbl{1}(\tilde{S}_{n-kd}\leq
r).
\]
We can show that the second term of the previous equation converges
to zero  a.s. Thus, we have $ \mathbbl{1}(\tilde{S}_{n}\leq
r)=\alpha_n$  a.s. Therefore,
\[
\tilde{S}_n=\mu_2+e_n+\sum_{i=1}^q\psi_ie_{n-i}+\Biggl[(\mu_1-\mu_2)+\sum_{i=1}^q(\phi_i-\psi_i)e_{n-i}\Biggr]\alpha_{n-d},
\qquad \mbox{a.s.},
\]
 that is,  $\tilde{S}_n=S_{n}$ a.s.  The proof is
complete.

%%%%%%%%%%%%%%%%%%%%%%%%%%%%%%%%%%%%%%%%%%%%%%%%%%%%%%%%%%%%%%%%%%%%%%%%%%%%%%%%%%%%%%%%%%%%%%%%%%%%%%%%%%%%%%%%%%%%%%%%%%
%                                                                                                                      %%%
%                                Appendix A.2: Proofs of Theorem 3.1                                                   %%%
%                                                                                                                      %%%
%%%%%%%%%%%%%%%%%%%%%%%%%%%%%%%%%%%%%%%%%%%%%%%%%%%%%%%%%%%%%%%%%%%%%%%%%%%%%%%%%%%%%%%%%%%%%%%%%%%%%%%%%%%%%%%%%%%%%%%%%%
%s5.2 #&#
\subsection{\texorpdfstring{Proof of Theorem \protect\ref{sethm2}}{Proof of Theorem 3.1}}
The notations $a_n$ and $b_n$ are defined by (\ref{op}), $m$ and
$\delta$ are the same as those in the proof of Theorem \ref{sethm1}.
From Theorem \ref{sethm1}, we have $y_n=a_n+(b_n-a_n)\alpha_{n-d}$.
For $n\geq (m+2)d+q$, we decompose $\alpha_{n-d}$ into two parts
\begin{eqnarray*}
\alpha_{n-d}&=&\sum_{j=1}^{[(n-q)/d]-1}\Biggl[\Biggl(\prod_{s=1}^{j-1}W_{n-sd}\Biggr)U_{n-jd}\Biggr]
+\sum_{j=[(n-q)/d]}^{\infty}\Biggl[\Biggl(\prod_{s=1}^{j-1}W_{n-sd}\Biggr)U_{n-jd}\Biggr]\\
&\equiv&I_1+I_2.
\end{eqnarray*}
Clearly, $I_1\in\mathcal{F}_d^{n-d}$ and
$I_2\in\mathcal{F}_{-\infty}^{n-d}$, where
$\mathcal{F}_m^n=\sigma(e_m,\ldots,e_n)$. By calculation, we have for
$n\geq (m+2)d+q$
\begin{eqnarray*}
|\operatorname{Cov}(y_0, y_n)|&\leq&
\sum_{j=[(n-q)/d]}^{\infty}\Biggl[\mathbb{E}\Biggl(\prod_{s=1}^{j-1}|W_{n-sd}|\Biggr)\Biggr]^{1/2}[\mathbb{E}(b_n-a_n)^2
(y_0-\mathbb{E}y_0)^2]^{1/2}\\
&\leq&\sum_{j=[(n-q)/d]}^{\infty}\sqrt{\delta}^{[(j-1)/(m+1)]}[\mathbb{E}(b_n-a_n)^2
]^{1/2}[\mathbb{E}(y_0-\mathbb{E}y_0)^2]^{1/2}\\
&\leq&\frac{H(m+1)}{1-\sqrt{\delta}}\sqrt{\delta}^{[(n-q-d)/(d(m+1))]-1}
\end{eqnarray*}
by H\"{o}lder's inequality, the boundedness of $W_n$ and $U_n$, and
the independence of $\{b_n, a_n\}$ and $y_0$, where
\[
H=\Biggl[|\mu_1-\mu_2|+(\mathbb{E}e_1^2)^{1/2}\sum_{i=1}^q|\phi_i-\psi_i|\Biggr]
\Biggl[|\mu_1|+|\mu_2|+(\mathbb{E}e_1^2)^{1/2}\sum_{i=1}^q(|\phi_i|+|\psi_i|)\Biggr].
\]
Thus, the conclusion holds.

%s5.3 #&#
\subsection{\texorpdfstring{Proof of Theorem \protect\ref{sethm3}}{Proof of Theorem 3.2}}
 Let $x_n=a_n+(b_n-a_n)I_1$. Then
$y_n-x_n=(b_n-a_n)I_2$. Clearly, $x_n\in \mathcal{F}_d^n$ and
$\mathbb{E}|y_n-x_n|=\mathrm{O}(\rho^n)$ for large enough $n$, where
$\rho\in(0, 1)$. So, using the independence of $x_n$ and $y_0$, for
large enough~$n$,  we have
\begin{eqnarray*}
&&|\mathbb{P}(y_0\leq u, y_n\leq
v)-\mathbb{P}(y_0\leq u)\mathbb{P}(y_n\leq v)|\\
&&\quad=|\mathbb{E}\{[\mathbbl{1}(y_n\leq v)-\mathbbl{1}(x_n\leq
v)][\mathbbl{1}(y_0\leq
u)-\mathbb{E}\mathbbl{1}(y_0\leq u)]\}|\\
&&\quad\leq\mathbb{E}|\mathbbl{1}(y_n\leq v)-\mathbbl{1}(x_n\leq v)|.
\end{eqnarray*}
On noting the independence between $e_n$ and $\bar{e}_{n-1}$, where
$
\bar{e}_{n-1}=\mu_2+\sum_{i=1}^q\psi_ie_{n-i}+(b_n-a_n)\alpha_{n-d}$,
the density of $y_n$ is $f_y(x)=\int_{\mathbb{R}}h(x-y)\,\mathrm{d}G_{\bar{e}}(y)$,
where $h(x)$ is the density function of $e_1$ and
$G_{\bar{e}}(y)$ is the distribution function of $\bar{e}_{n-1}$. On
using the property of convolution, $f_y(x)$ is continuous and
bounded. Write $\|f_y\|_{\infty}=\max\{|f_y(x)|\dvt x\in\mathbb{R}\}$.
On the one hand, using the following inequality
\[
|\mathbbl{1}(x\leq t)-\mathbbl{1}(y\leq t)|\mathbbl{1}(|x-y|\leq
\epsilon)\leq \mathbbl{1}(t-\epsilon\leq x\leq t+\epsilon),
\]
we can get
\[
\mathbb{E}[|\mathbbl{1}(y_n\leq v)-\mathbbl{1}(x_n\leq
v)|\mathbbl{1}(|y_n-x_n|\leq \epsilon)]\leq
\mathbb{P}(v-\epsilon\leq y_n\leq v+\epsilon)
\leq2\|f_y\|_\infty\epsilon.
\]
On the other hand, using Markov's inequality, we have
\[
\mathbb{E}\{|\mathbbl{1}(y_n\leq v)-\mathbbl{1}(x_n\leq
v)|\mathbbl{1}(|y_n-x_n|>\epsilon)\}\leq
\mathbb{E}\mathbbl{1}(|y_n-x_n|>\epsilon)\} \leq\epsilon^{-1}
\mathbb{E}|y_n-x_n|.
\]
Choosing $\epsilon=\mathrm{O}(\rho^{n/2})$, we can obtain
\[
|\mathbb{P}(y_0\leq u, y_n\leq v)-\mathbb{P}(y_0\leq
u)\mathbb{P}(y_n\leq v)|=\mathrm{O}(\rho^{n/2}).
\]
Hence, the result holds.
\end{appendix}

\section*{Acknowledgements}
We thank the referee, the Associate Editor and the Editor for their
very helpful comments and suggestions.  The research was partially
supported by Hong Kong Research Grants Commission Grants HKUST601607
and HKUST602609, the Saw Swee Hock professorship at the National University of
Singapore and the Distinguished Visiting Professorship at the University of
Hong Kong (both to HT).

% imsref loaded by svajune.rapalyte, 2011-10-04 11:10:26
% imsref loaded by svajune.rapalyte, 2011-10-05 10:14:45

\printhistory


\begin{thebibliography}{19}
% BibTex style file: bej.bst, 2010-01-21
% Default style options (sort=1,type=number).
% Used options (sort=1,type=number).

%b1 #&#
\bibitem{ref1}
\begin{barticle}[mr]
\bauthor{\bsnm{An},~\bfnm{H.~Z.}\binits{H.Z.}} \AND
  \bauthor{\bsnm{Huang},~\bfnm{F.~C.}\binits{F.C.}}
(\byear{1996}).
\btitle{The geometrical ergodicity of nonlinear autoregressive models}.
\bjournal{Statist. Sinica}
\bvolume{6}
\bpages{943--956}.
\bid{issn={1017-0405}, mr={1422412}}
\end{barticle}
\endbibitem

%b2 ###
%b2 #&#
\bibitem{ref2}
\begin{barticle}[mr]
\bauthor{\bsnm{Brockwell},~\bfnm{Peter~J.}\binits{P.J.}},
  \bauthor{\bsnm{Liu},~\bfnm{Jian}\binits{J.}} \AND
  \bauthor{\bsnm{Tweedie},~\bfnm{Richard~L.}\binits{R.L.}}
(\byear{1992}).
\btitle{On the existence of stationary threshold autoregressive moving-average
  processes}.
\bjournal{J. Time Ser. Anal.}
\bvolume{13}
\bpages{95--107}.
\bid{doi={10.1111/j.1467-9892.1992.tb00096.x}, issn={0143-9782}, mr={1165659}}
\end{barticle}
\endbibitem

%b3 ###
%b3 #&#
\bibitem{ref3}
\begin{barticle}[mr]
\bauthor{\bsnm{Chan},~\bfnm{K.~S.}\binits{K.S.}},
  \bauthor{\bsnm{Petruccelli},~\bfnm{Joseph~D.}\binits{J.D.}},
  \bauthor{\bsnm{Tong},~\bfnm{H.}\binits{H.}} \AND
  \bauthor{\bsnm{Woolford},~\bfnm{Samuel~W.}\binits{S.W.}}
(\byear{1985}).
\btitle{A multiple-threshold {AR{$(1)$}} model}.
\bjournal{J. Appl. Probab.}
\bvolume{22}
\bpages{267--279}.
\bid{issn={0021-9002}, mr={0789351}}
\end{barticle}
\endbibitem

%b4 ###
%b4 #&#
\bibitem{ref4}
\begin{barticle}[mr]
\bauthor{\bsnm{Chan},~\bfnm{K.~S.}\binits{K.S.}} \AND
  \bauthor{\bsnm{Tong},~\bfnm{H.}\binits{H.}}
(\byear{1985}).
\btitle{On the use of the deterministic {L}yapunov function for the ergodicity
  of stochastic difference equations}.
\bjournal{Adv. in Appl. Probab.}
\bvolume{17}
\bpages{666--678}.
\bid{doi={10.2307/1427125}, issn={0001-8678}, mr={0798881}}
\end{barticle}
\endbibitem

%b5 ###
%b5 #&#
\bibitem{ref5}
\begin{barticle}[mr]
\bauthor{\bsnm{Chen},~\bfnm{Rong}\binits{R.}} \AND
  \bauthor{\bsnm{Tsay},~\bfnm{Ruey~S.}\binits{R.S.}}
(\byear{1991}).
\btitle{On the ergodicity of {${\rm TAR}(1)$} processes}.
\bjournal{Ann. Appl. Probab.}
\bvolume{1}
\bpages{613--634}.
\bid{issn={1050-5164}, mr={1129777}}
\end{barticle}
\endbibitem

%b6 ###
%b6 #&#
\bibitem{ref6}
\begin{barticle}[mr]
\bauthor{\bsnm{Cline},~\bfnm{Daren B.~H.}\binits{D.B.H.}} \AND
  \bauthor{\bsnm{Pu},~\bfnm{Huay-min~H.}\binits{H.m.H.}}
(\byear{1999}).
\btitle{Geometric ergodicity of nonlinear time series}.
\bjournal{Statist. Sinica}
\bvolume{9}
\bpages{1103--1118}.
\bid{issn={1017-0405}, mr={1744827}}
\end{barticle}
\endbibitem

%b7 ###
%b7 #&#
\bibitem{ref7}
\begin{barticle}[mr]
\bauthor{\bsnm{Cline},~\bfnm{Daren B.~H.}\binits{D.B.H.}} \AND
  \bauthor{\bsnm{Pu},~\bfnm{Huay-Min~H.}\binits{H.M.H.}}
(\byear{2004}).
\btitle{Stability and the {L}yapounov exponent of threshold {AR}-{ARCH}
  models}.
\bjournal{Ann. Appl. Probab.}
\bvolume{14}
\bpages{1920--1949}.
\bid{doi={10.1214/105051604000000431}, issn={1050-5164}, mr={2099657}}
\end{barticle}\vadjust{\goodbreak}
\endbibitem



%b9 ###
%b9 #&#
\bibitem{ref9}
\begin{barticle}[mr]
\bauthor{\bsnm{Ling},~\bfnm{Shiqing}\binits{S.}}
(\byear{1999}).
\btitle{On the probabilistic properties of a double threshold {ARMA}
  conditional heteroskedastic model}.
\bjournal{J. Appl. Probab.}
\bvolume{36}
\bpages{688--705}.
\bid{issn={0021-9002}, mr={1737046}}
\end{barticle}
\endbibitem


%b11 ###
%b11 #&#
\bibitem{ref11}
\begin{barticle}[mr]
\bauthor{\bsnm{Ling},~\bfnm{Shiqing}\binits{S.}},
  \bauthor{\bsnm{Tong},~\bfnm{Howell}\binits{H.}} \AND
  \bauthor{\bsnm{Li},~\bfnm{Dong}\binits{D.}}
(\byear{2007}).
\btitle{Ergodicity and invertibility of threshold moving-average models}.
\bjournal{Bernoulli}
\bvolume{13}
\bpages{161--168}.
\bid{doi={10.3150/07-BEJ5147}, issn={1350-7265}, mr={2307400}}
\end{barticle}
\endbibitem

%b12 ###
%b12 #&#
\bibitem{ref12}
\begin{barticle}[mr]
\bauthor{\bsnm{Liu},~\bfnm{Jian}\binits{J.}} \AND
  \bauthor{\bsnm{Susko},~\bfnm{Ed}\binits{E.}}
(\byear{1992}).
\btitle{On strict stationarity and ergodicity of a nonlinear {ARMA} model}.
\bjournal{J.~Appl. Probab.}
\bvolume{29}
\bpages{363--373}.
\bid{issn={0021-9002}, mr={1165221}}
\end{barticle}
\endbibitem

%b13 ###
%b13 #&#
\bibitem{ref13}
\begin{barticle}[mr]
\bauthor{\bsnm{Lu},~\bfnm{Zudi}\binits{Z.}}
(\byear{1998}).
\btitle{On the geometric ergodicity of a non-linear autoregressive model with
  an autoregressive conditional heteroscedastic term}.
\bjournal{Statist. Sinica}
\bvolume{8}
\bpages{1205--1217}.
\bid{issn={1017-0405}, mr={1666249}}
\end{barticle}
\endbibitem

%b14 ###
%b14 #&#
\bibitem{ref14}
\begin{barticle}[mr]
\bauthor{\bsnm{Robinson},~\bfnm{P.~M.}\binits{P.M.}}
(\byear{1977}).
\btitle{The estimation of a nonlinear moving average model}.
\bjournal{Stochastic Process. Appl.}
\bvolume{5}
\bpages{81--90}.
\bid{issn={0304-4149}, mr={0428654}}
\end{barticle}
\endbibitem

%b15 ###
%b15 #&#
\bibitem{ref15}
\begin{bmisc}[auto:STB|2011/09/12|07:03:23]
\bauthor{\bsnm{Slutsky},~\bfnm{E.}\binits{E.}}
(\byear{1927}).
\bhowpublished{The summation of random causes as the source of cyclic
  processes. \textit{Voprosy Koniunktury} \textbf{3} 34--64. \{English
  translation in \textit{Econometrika} \textbf{5} 105--146.\}}.
\end{bmisc}
\endbibitem

%b16 ###
%b16 #&#
\bibitem{ref16}
\begin{bincollection}[auto:STB|2011/09/12|07:03:23]
\bauthor{\bsnm{Tong},~\bfnm{Howell}\binits{H.}}
(\byear{1978}).
\btitle{On a thresold model}.
In \bbooktitle{Pattern Recognition and Signal Processing}
(\beditor{\bfnm{Chi Hau}\binits{C.H.}~\bsnm{Chen}}, ed.)
\bpages{575--586}.
\baddress{Amsterdam}: \bpublisher{Sijthoff and Noordhoff}.
\end{bincollection}
\endbibitem


%b17 ###
%b17 #&#
\bibitem{ref17}
\begin{bbook}[mr]
\bauthor{\bsnm{Tong},~\bfnm{Howell}\binits{H.}}
(\byear{1990}).
\btitle{Nonlinear Time Series: A Dynamical System Approach}.
\bseries{Oxford Statistical Science Series}
\bvolume{6}.
\baddress{New York}: \bpublisher{Oxford Univ. Press}.
\bid{mr={1079320}}
\end{bbook}
\endbibitem

%b18 ###
%b18 #&#
\bibitem{ref18}
\begin{bmisc}[auto:STB|2011/09/12|07:03:23]
\bauthor{\bsnm{Tong},~\bfnm{H.}\binits{H.}}
(\byear{2011}).
\bhowpublished{Threshold models in time series analysis -- 30 years on.
  \textit{Stat. Interface} \textbf{4} 107--118}.
\end{bmisc}
\endbibitem

\end{thebibliography}
\end{document}